\documentclass[twocolumn]{autart}    % Enable this line and disable the 

\usepackage{graphicx}  % 用于插入图片
\usepackage{subcaption} % 用于子图                                 % preceding line to obtain a two-column 
\usepackage{cite}
\usepackage{amsmath,amssymb,amsfonts}
\usepackage{algorithmic}
\usepackage{graphicx}
\usepackage{amsmath,amssymb}
\usepackage{amsfonts}
\usepackage{lmodern}
\usepackage{textcomp}
\usepackage{xcolor}                                  
\usepackage{graphicx} 
                              \newtheorem{theorem}{\bf{Theorem}}

                              \newtheorem{assumption}{Assumption}

                              \newtheorem{lemma}{\bf{Lemma}}

                              \newenvironment{proof}{\noindent{\em Proof:\/}}{\hfill $\Box$\par}
                              \newtheorem{remark}{\bf{Remark}}

\begin{document}

\begin{frontmatter}

\title{ Distributed Stochastic Optimization for Non-Smooth and Weakly Convex Problems under Heavy-Tailed Noise\thanksref{footnoteinfo}} % Title, preferably not more 

\thanks[footnoteinfo]{This work was supported in part by the Joint Funds of the National
Natural Science Foundation of China under Grant U24A20258, in part by the National Natural Science Foundation of China under Grants 92367205, in part by the Zhejiang Provincial Natural Science Foundation of China under Grant LRG25F030001, in part by the funding of Leading Innovative and Entrepreneur
Team Introduction Program of Zhejiang under Grant 2023R01006, in part by the Fundamental Research Funds for the Provincial Universities of Zhejiang under Grant RF-C2023007.\\ }
\thanks{Corresponding author: Bo Chen (bchen@aliyun.com).}
\author[a,b]{Jun Hu},
\author[a,b]{Chao Sun},
\author[a,b]{Bo Chen}, 
\author[a,b]{Jianzheng Wang},  % Add the 
\author[a,b]{Zheming Wang} 

\address[a]{Department of Automation, Zhejiang University of Technology, Hangzhou 310023, China}  % Please supply                                                    % here.
\address[b]{Zhejiang Key Laboratory of Intelligent Perception and Control for Complex Systems, Hangzhou 310023, China}
          
                             % keyword list or with the 
                                          % help of the Automatica 
                                          % keyword wizard

\begin{abstract}                          % Abstract of not more than 200 words.
    In existing distributed stochastic optimization studies, it is usually assumed that the gradient noise has a bounded variance. However, recent research shows that the heavy-tailed noise, which allows an unbounded variance, is closer to practical scenarios in many tasks. Under heavy-tailed noise, traditional optimization methods, such as stochastic gradient descent, may have poor performance and even diverge. Thus, it is of great importance to study distributed stochastic optimization algorithms applicable to the heavy-tailed noise scenario. However, most of the existing distributed algorithms under heavy-tailed noise are developed for convex and smooth problems, which limits their applications. This paper proposes a clipping-based distributed stochastic algorithm under heavy-tailed noise that is suitable for non-smooth and weakly convex problems. The convergence of the proposed algorithm is proven, and the conditions on the parameters are given. A numerical experiment is conducted to demonstrate the effectiveness of the proposed algorithm.
\end{abstract}

\end{frontmatter}

\section{INTRODUCTION}
\hspace{1em}In recent years, distributed optimization has gained significant attention due to its wide applications in large-scale machine learning\cite{verbraeken2020survey}, robotic networks\cite{mueller2019modern}, and energy systems\cite{chen2017delay}. Unlike centralized systems, 
a number of interconnected agents cooperate in order to achieve a desirable global objective. Distributed optimization overcomes a single point of failure, high communication requirements, substantial computation burden, and limited flexibility and scalability\cite{yang2019survey}.\par
\hspace{1em}Stochastic optimization refers to solving the problem in the form of the expectation of a stochastic function. The most popular algorithm is stochastic gradient descent (SGD)\cite{amiri2020machine,meng2019convergence,shamir2014distributed,george2020distributed,swenson2022distributed,pu2021sharp,chen2021distributted}. For certain applications in large-scale distributed machine learning, the authors in \cite{yu2019linear} presented a distributed parallel SGD with momentum for training transformer networks\cite{vaswani2017attention}. \cite{yi2022primal} proposed a distributed primal-dual SGD algorithm for neural networks. \cite{assran2019stochastic} introduced stochastic gradient push which enables the use of generic communication topologies for training ResNet.\par
\hspace{1em}To utilize SGD algorithm in large-scale optimization, a critical challenge in distributed SGD arises from the inherent noise in gradient estimation for each agent \cite{dhasunr2010distributed}. Existing works usually assume that the variance of the noise, which refers to the difference between the estimated stochastic gradient and the true gradient is bounded\cite{yu2019linear,yi2022primal,assran2019stochastic}. However, the finite variance assumption seems too optimistic, and the heavy tail assumption is more realistic. For example, the authors in \cite{gurbuzbalaban2021heavy} pointed out that the iterates appear heavy-tail behavior in SGD algorithm, which means the noise can be non-Gaussian. \cite{armacki2023high,garg2021proximal,battash2024revisiting} showed that the gradient noise follows a heavy-tailed distribution which violate the bounded variance assumption. 
\cite{gurbuzbalaban2024heavy} demonstrated the heavy-tail behavior of decentralized stochastic gradient descent and get more precise characterizations of the tail behavior. The authors in \cite{zhang2020adaptive} observed the noise distribution when training BERT\cite{bao2021beit} on the Wikipedia and noticed that the noise follow Lévy-$\alpha$-stable distribution, which is a heavy-tailed distribution. Moreover, in the presence of heavy-tailed noise, the SGD algorithm frequently exhibits poor performance and even diverges.\par
\hspace{1em}To address the poor performance of SGD under heavy-tailed noise, gradient clipping-based algorithms have been studied in many literature. For example, \cite{zhang2020adaptive} proposed a gradient clipping-based method for strongly convex and weakly convex optimization, and showed the convergence in expectation. In \cite{li2022high}, the authors showed high probability convergence for stochastic optimization under heavy-tailed noise. For distributed cases, in \cite{yang2022taming}, the authors studied stochastic optimization algorithms for federated learning under heavy-tailed noise.
In \cite{liu2022communication}, the authors designed a communication-efficient local gradient clipping algorithm.\par
\hspace{1em}The weakly convex problem\cite{davis2019stochastic}, which is a special type of non-convex problem, has many important applications in  machine learning, such as phase retrieval \cite{duchi2019solving}, Generative Adversarial Nets (GANs)\cite{liu2021first}, and  input weakly convex neural networks (IWCNN)\cite{shumaylov2024weakly}. The authors in \cite{gao2023delayed} proposed a  delayed stochastic algorithm for weakly convex optimization via a network connected to a master node. In addition,
the authors in \cite{chen2021distributed} showed a linear convergence rate under a sharpness condition for the weakly convex problem. However, all of these studies are developed for deterministic problems or assume that the stochastic noise satisfies the bounded variance condition.  \par
\hspace{1em}In our work, we investigate a distributed stochastic optimization method for non-smooth and weakly convex problems affected by heavy-tailed noise. The key contributions of this paper include:\par
(1)We design a distributed stochastic subgradient descent algorithm based on gradient clipping method to solve the weakly convex and non-smooth problem under heavy-tailed noise. We present a rigorous proof for convergence analysis and also show the convergence rate under some mild parameter conditions.\par
(2)In contrast to \cite{sun2023distributed,yang2024online}, which study distributed stochastic optimization problems under heavy-tailed noises, our algorithm can be applied to  non-smooth and weakly convex optimization equipped with weaker assumptions.\par
(3)Compared with the existing works on  stochastic optimization in heavy-tail phenomenon\cite{gorbunov2023high,yang2022taming,liu2022communication}, 
our approach is distributed in the sense that agents communicate only with neighboring nodes via a communication graph, thereby obviating the requirement for a centralized coordination server. While \cite{gorbunov2023high,yang2022taming} depend on central server, \cite{liu2022communication} involves averaging agents' states at specific time steps, 
our work explicitly addresses multi-agent consensus based on graph theory. \\

\section{NOTATIONS AND PRELIMINARIES}\label{section2}
\hspace{1em}Throughout this paper, $\mathbb{R}^n$ denotes the set of all real-valued vectors of dimension $n$. The ``$\inf$" represents the greatest lower bound (infimum) of a given sequence. $||\cdot||$ indicates the standard Euclidean norm (2-norm) of a vector. $\mathcal{P}_\Omega$ denotes the orthogonal projection operator onto the set 
$\Omega$. The notation $[A(k)]_{i,j}$ denotes the $(i,j)$-th entry of the adjacency matrix at iteration $k$. $\mathbb{E}[\cdot]$ represents the expectation of a random variable. $\partial f(x) $ denotes the subgradient of $f$ at point $x$.
\subsection{Network}
\hspace{1em} In distributed optimization, agents communicate over a time-varying network represented by a graph $\mathcal{G}(k) = (\mathcal{V},\mathcal{E}_k)$, where $\mathcal{V}=\{1,\dots,N\}$ is the set of nodes, and $\mathcal{E}_k\subseteq\mathcal{V}\times\mathcal{V}$ denotes the edge set at each iteration $k>0$. Let $A(k)=[a_{i,j}(k)]$ denote  weighted adjacency matrix  associated with $\mathcal{G}(k)$. For each node $i$, the neighborhood at time $k$ is defined as $\mathcal{N}_i(k)=\{j\mid a_{i,j}(k)>0\}$.
\subsection{Weak Convexity}
\hspace{1em} A function $f:\mathbb{R}^N\to\mathbb{R}\cup\{+\infty\}$ is said to be $\rho$-weakly convex, if the assignment $x \mapsto f(x) + \frac{\rho}{2} \|x\|^2 $
is a convex function. We define the subgradient of $f(x)$ as\par
\begin{equation}
    \partial f(x)= \partial h(x)- \rho||x||, \label{0}
\end{equation}
where $h(x)$ is a convex function.\par
\hspace{1em} According to the definition of weak convexity, the subgradient of $f$ at $x$, denoted by $\partial f(x)\subset \mathbb{R}^N$, satisfies
\begin{equation}
f(y)\ge f(x)+\langle \partial f(x),y-x \rangle-\frac{\rho}{2}||y-x||^2. \label{1}
\end{equation}
\hspace{1em}The following lemma provides a useful property of weak convexity, which closely resembles the definition of convexity and plays a key role in our convergence analysis.\par
\begin{lemma}
    \cite[Lemma II.1]{chen2021distributed}: Assume that \( f(x) \) is \( \rho \)-weakly convex in \( \mathbb{R}^n \). Then for \( \forall x_1, \ldots, x_k \in \mathbb{R}^n \), it follows that
\begin{equation}
f\left( \sum_{i=1}^k \omega_i x_i \right) \leq \sum_{i=1}^k \omega_i f(x_i) + \frac{\rho}{2} \sum_{1\le i \le j \le k} \omega_i \omega_j \| x_i - x_j \|^2,
\end{equation}
where $\sum_{i=1}^{k}\omega_i=1$ and $ \omega_i \ge 0 \text{ for }\forall i$.\label{lemma1}
\end{lemma}
\subsection{Moreau Envelope and Proximal Mapping}
\hspace{1em}To study the stochastic subgradient method, the principal goal of non-smooth and weakly convex optimization 
is to find a stationary point, \emph{i.e.,} a point $x\in\Omega$ that satisfies $0\in\partial f(x)$. An effective method to evaluate the stationarity metric is to construct the Moreau envelope.  Let's define the standard form \cite{davis2019stochastic}\cite{klee1971convex} for a $\rho$-weakly convex function $f$ with $\varphi(x)=f(x)+\mathbb{I}_{\Omega}(x)$
\footnote{In order to restrict $dom f$ to set $\Omega$, $\mathbb{I}_{\Omega}(x)$ satisfies $\mathbb{I}_{\Omega}(x)=0$ when $x\in\Omega$, and $\mathbb{I}_{\Omega}(x)=\infty$ otherwise.}
\begin{equation}
\varphi_{\mu}(x):=\min_{y \in \Omega} \varphi(y)+\frac{1}{2\mu}||y-x||^2,1<\mu<1/\rho. \label{3}
\end{equation}
\hspace{1em}Suppose the function $f$ satisfies the Lipschitz continuity condition with parameter $L>0$ that
\begin{equation}
    \| f(x) - f(y) \| \leq L \| x - y \|\quad x,y\in\Omega.\label{4}
\end{equation}   

According to \cite[Lemma 4.2]{drusvyatskiy2019efficiency}, $\varphi(y)+\frac{1}{2\mu}||y-x||^2$ is strictly convex,
while the Moreau envelope is a $C^1$ smoothed version (first-order continuously differentiable) of the non-smooth function $f(x)$ over $\Omega$. Define 
\begin{equation}
\hat{x}=\mathop{\text{argmin}}\limits_{y \in \Omega}\varphi(y)+\frac{1}{2\mu}||y-x||^2. \label{5}
\end{equation}
The Moreau envelope $\varphi_{\mu}$ has gradient
\begin{equation}
\nabla\varphi_{\mu}(x)=\frac{1}{\mu}(x-prox_{\mu f}(x)),\label{6}
\end{equation}
where $prox_{\mu f}(x):=\hat{x}$ is the proximal mapping \cite{chen2021distributed}. $prox_{\mu f}$ is well-defined and single-valued\cite[Lemma 4.3]{drusvyatskiy2019efficiency}.
The proof of \eqref{6} can be found in Appendix \ref{AppendixA}.\par
\hspace{1em}Based on the optimality condition, we have
\begin{equation}
0\in \partial\varphi(\hat{x})+\frac{1}{\mu}||\hat{x}-x||.\label{7}
\end{equation}
Furthermore, the following stationarity metric holds
\begin{equation}
dist(0,\partial\varphi(\hat{x}))\le \frac{1}{\mu}||\hat{x}-x||.\label{8}
\end{equation}
Taking the norm on both sides of \eqref{6}, we have that a small $||\nabla\varphi_{\mu}(x)||$ implies that $x$ is near the point $\hat{x}$. Based on \eqref{8}, when $||\hat{x}-x||$ is small, $dist(0,\partial\varphi(\hat{x}))$ is small. According to the optimality condition, if $dist(0,\partial\varphi(\hat{x}))=0$, then $\hat{x}$ is a  stationary point of $f(x)$. Thus, a small $||\nabla\varphi_{\mu}(x)||$ implies that $x$ is close to a stationary point of $f(x)$.  In the following analysis, we will demonstrate the decay of $||\nabla\varphi_{\mu}(x)||$ along the sequence generated by our algorithm.\par
\hspace{1em} The following lemma shows a property of the proximal mapping.\par
\begin{lemma}\cite[Lemma II.8]{chen2021distributed}(Proximal Mapping):Assume that $f(x)$ is $\rho$-weakly convex. Then, there exists $\mu<1/\rho$ satisfies
\begin{equation}
    \begin{split}
       ||prox_{\mu f}(x_1)-prox_{\mu f}(x_2)||\le \frac{1}{1-\mu\rho}||x_1-&x_2||\ \\
                          &\forall x_1,x_2\in\Omega.
    \end{split}
\end{equation}\label{lemma2}
\end{lemma}
\section{PROBLEM FORMULATION}
\hspace{1em}In this section, we  introduce the problem description of distributed stochastic optimization.\par
\hspace{1em}Consider a multi-agent network consisting of $N>1$ agents collaborating to address the following problem:
\begin{equation}
\min\limits_{{x} \in \Omega} f(x) = \frac{1}{N} \sum_{i=1}^N \underbrace{\mathbb{E}_{\xi_i \sim \mathcal{D}_i} \left[ F_i(x; \xi_i) \right]}_{\triangleq f_i(x)} \label{10}
\end{equation}
where $x\in \Omega \subset \mathbb{R}^n$, $\Omega$ is the local constraint set,
$f(x)$ is the global objective function, $f_i(x)\triangleq\mathbb{E}_{\xi_i\sim D_i}[F_i(x,\xi_i)]$ is the local objective function for agent $i$,
and each $\xi_i$ following distribution $\mathcal{D}_i$ represents the local dataset for agent $i$. Moreover, $f_i(x)$ has subgradient $\partial f_i(x)$.\par
\hspace{1em}
At each iteration, every agent computes a stochastic gradient based on its local information. In other words, each agent has access to a stochastic first-order oracle$(\mathcal{SFO})$\cite{xin2019distributed},
which offers a stochastic subgradient $g_{i,k}(\cdot,\xi_{i,k})$ at each iteration, where ``$\cdot$" is the base point for agent $i$.\par
\hspace{1em}The following assumptions are widely adopted in the literature on distributed stochastic optimization(\emph{e.g.,}\cite{chen2021distributed},\cite{nedic2009distributed}).\par
\begin{assumption}
(Communication Network):The graph $G_N$ is strongly connected. The adjacency matrix $A$ is doubly stochastic, \textit{i.e.,} 
 $\sum_{i=1}^N a_{i,j}(k) = \sum_{j=1}^N a_{i,j}(k) = 1$. There exists a scalar $0 < \eta < 1$ 
 such that $a_{i,j}(k) \geq \eta$ when $j \in \mathcal{N}_i$. \label{assumption1}
\end{assumption}
\begin{assumption}
The graph $(\mathcal{V},E_{\infty})$ is strongly connected, where $E_{\infty}$ is the set of edges $(j,i)$ representing agent pairs communicating
 directly infinitely many times.\label{assumption2}
\end{assumption}
\begin{assumption}
    {(Bounded Intercommunication Interval):} There exists an integer $B \geq 1$ such that for every $(j, i) \in E_\infty$, agent $j$ sends its information to the neighboring agent $i$ at least once every $B$ consecutive time slots, \textit{i.e.,} at time $t_k$ or at time $t_k + 1$ or \dots or (at latest) at time $t_k + B - 1$ for any $k \geq 0$.\label{assumption3}
\end{assumption}
\begin{assumption}
The constraint set $\Omega$ is compact, convex and nonempty.\label{assumption4}
\end{assumption}
\begin{remark}
    Assumption \ref{assumption4} is commonly used in distributed optimization algorithms, \textit{e.g.,} \cite{chen2021distributed},\cite{liang2019distributed}. It means that the decision variable should be limited within a certain bounded region.
\end{remark}
\begin{assumption}
Suppose each local function $f_i$ is  a non-smooth and $\rho$-weakly convex function. By condition \eqref{1}, the global function $f(x)$ is also $\rho$-weakly convex. \label{assumption5}
\end{assumption}
\begin{assumption}
Suppose each local function $f_i$ satisfies $L$-lipschitz continuity as shown in \eqref{4}.\label{assumption6}
\end{assumption}
\begin{remark}
     In contrast to \cite{yang2024online}, which relies on the standard 
     $L$-smoothness assumption for analyzing non-convex problems, our work lifts this restriction by addressing non-smooth settings, where the objective functions may not have Lipschitz-continuous gradient. Moreover, different from \cite{sun2023distributed}, which focuses on distributed stochastic algorithms for smooth convex functions, we consider a more general setting of non-smooth and weakly convex objectives under a distributed subgradient framework.
\end{remark}
\section{ALGORITHM DESIGN}
\hspace{1em}In this section, we propose the algorithm, and prove its convergence. That is, if  Assumptions \ref{assumption1}-\ref{assumption8} hold, the sequence $\{ \varphi_{\mu}(\bar{x}_k)\}$ converges and  $ \mathbb{E}[\inf_{k_0\le k\le T}||\nabla\varphi_{\mu}(\bar{x}_k)||^2]$ converges to zero. 
Finally, we provide the convergence rate under certain conditions.
\subsection{Algorithm}
\hspace{1em}Motivated by \cite{zhang2020adaptive},\cite{sun2023distributed}, the updating law for the $i$-th agent is designed as follows
\begin{equation}
x_{i,k+1} =   \mathcal{P}_{\Omega}[v_{i,k} - \alpha_k \hat{g}_{i,k}], \label{14}
\end{equation}
where 
\[
v_{i,k} = \sum_{j=1}^N a_{i,j}(k) x_{j,k},
\]
and
\begin{align*}
    \hat{g}_{i,k} &= \min \left\{ 1, \frac{\tau_k}{\| g_{i,k}(v_{i,k},\xi_{i,k}) \|} \right\} g_{i,k}(v_{i,k},\xi_{i,k}).\\
\end{align*}
Here, \( x_{i,k} \) denotes the local decision variable of agent \( i \) at iteration \( k \), and \( v_{i,k} \) is an intermediate variable obtained by taking a weighted average of the neighboring \( x_{j,k} \).  
The sequences \( \alpha_k \) and \( \tau_k \) are positive and will be specified later.  
In addition, \( \xi_{i,k} \) represents a random variable corresponding to the mini-batch sampled by agent \( i \) at iteration \( k \). \par
\hspace{1em}Let \(\mathcal{F}_k\) denote the \(\sigma\)-algebra generated by the sequence of random variables from \(0\) to \(k-1\), \emph{i.e.,} $\mathcal{F}_k=\{\xi_{i,m}|i\in\mathcal{V},m=0,1,\dots k-1\}$.\par
\hspace{1em}To facilitate the subsequent analysis, we adopt the following assumptions about the noise.\par
\begin{assumption}
(Unbiased Local Gradient Estimator):
$\mathbb{E}[g_{i,k}(v_{i,k},\xi_{i,k})|\mathcal{F}_k]\in\partial f_i(v_{i,k})$ with probability 1. Without loss of generality, let $\mathbb{E}[g_{i,k}|\mathcal{F}_k]=G_{i,k}\in \partial f_i(v_{i,k})$. \par
\hspace{1em}Let $b_{i,k}=\hat{g}_{i,k}-G_{i,k}$ and $B_{i,k}=E[b_{i,k}|\mathcal{F}_k]=
\mathbb{E}[\hat{g}_{i,k}|\mathcal{F}_k]-G_{i,k}$.\label{assumption7}
\end{assumption}
\begin{assumption} 
There exist two positive constants $\alpha \in (1,2] $
and $\gamma >0$ such that for $G_{i,k}\in \partial f_i(v_{i,k})$, $\mathbb{E}[||g_{i,k}(v_{i,k},\xi_{i,k})-G_{i,k}||^{\alpha}|\mathcal{F}_k]\le \gamma^{\alpha}$. The noise is 
said to be heavy-tailed if $\alpha<2 $.\label{assumption8}
\end{assumption}
\hspace{1em} According to the definition of the subgradient for $\rho$-weakly convex function in \eqref{0}, and given that every $f_i(x)$ is $L$-lipschitz on compact set $\Omega$, subgradient $||\partial f_i(x)||$ is uniformly bounded on $\Omega$\cite{rockafellar2009variational}, \emph{i.e.,}
$||\partial f_i(x)||\le C_0$ for $x\in \Omega$. According to the definition of $v_{i,k}$, Assumption \ref{assumption1} and Assumption \ref{assumption4}, $v_{i,k}\in \Omega$ for every iteration $k$. Thus, $||\partial f_i(v_{i,k})||\le C_0$ and $\Vert G_{i,k} \Vert\le C_0$ for all $k$. \par
\begin{remark}
In contrast to \cite{chen2021distributed}, our work explicitly considers the presence of heavy-tailed noise in the stochastic gradient. Inspired by the clipping technique introduced in  \cite{zhang2020adaptive}, we incorporate a gradient clipping mechanism into the design of a  distributed subgradient algorithm to address non-smooth and weakly convex optimization problems under such noise conditions.
\end{remark}
\subsection{Some Technical Lemmas}
\hspace{1em}Before the convergence analysis, we introduce some key lemmas.\par
\begin{lemma} Let Assumptions \ref{assumption7} and \ref{assumption8} hold. Then, for any agent $i$, if $\tau_k \ge 2C_0 $, 
we have
\begin{align*}
    &||B_{i,k}||=||\mathbb{E}[\hat{g}_{i,k}|\mathcal{F}_k]-G_{i,k}||\le (2\gamma)^\alpha\tau_k^{1-\alpha}. \label{15}
\end{align*}
\begin{proof}
See Appendix \ref{sectionB}.\label{lemma5}
\end{proof}
\end{lemma}
\hspace{1em}The following lemma illustrates the improvement after each iteration \eqref{14}.
\begin{lemma}
Under Assumptions \ref{assumption1}-\ref{assumption8}, if $\mu<\frac{1}{2\rho}$, we have
\begin{equation}
    \begin{split}
      &\sum_{i=1}^{N}||x_{i,k+1}-\hat{v}_{i,k}||\\
      \le&\quad\sum_{i=1}^{N}||v_{i,k}-\hat{v}_{i,k}||^2\\
      &\quad+2\alpha_k \bigg(N\left(-\frac{1}{2\mu}+\rho\right)||\bar{x}_k-z_k||^2\\
      &\quad+\frac{L(2-\mu\rho)}{1-\mu\rho}\sum_{i=1}^{N}||\bar{x}_k-x_{i,k}||\\
      &\quad+2\rho\left(1+\frac{1}{(1-\mu\rho)^2}\right)\sum_{i=1}^{N}||\bar{x}_k-x_{i,k}||^2\\
      &\quad-\sum_{i=1}^{N}b_{i,k}^{\top}(v_{i,k}-\hat{v}_{i,k})\bigg)+N\tau_k^2\alpha_k^2, \label{16}
    \end{split}
\end{equation} 
where $\bar{x}_k=1/N\sum_{i=1}^{N}x_{i,k}$, $z_k := \mathop{\text{argmin}}\limits_{y \in \Omega}  f(y) + \frac{1}{2\mu} \|y - \bar{x}_k\|^2 $,
$\hat{v}_{i,k}:= \mathop{\textrm{argmin}}\limits_{y \in \Omega} f(y) + \frac{1}{2\mu} \|y - v_{i,k}\|^2 $ 
.\par
\begin{proof}     
See Appendix \ref{sectionC}. \label{lemma6}
\end{proof}
\end{lemma}
\hspace{1em}Before establishing the convergence of the global objective, it is essential to show that the consensus error goes to zero.
\begin{lemma} 
Under Assumptions \ref{assumption1}-\ref{assumption3}, for the iteration \eqref{14}, if the stepsize satisfies $lim_{k\to\infty}\alpha_{k+1}/\alpha_k=1$, then
\[
\lim_{k\to\infty}||\bar{x}_k-x_{i,k}||=0\quad \forall i.
\]
Denote $\Delta_k\in R^{Nn}$ as a vector, where the $i$-th element $\Delta_{k,i}=||\bar{x}_k-x_{i,k}||$. Then, we have
\begin{equation}
    ||\Delta_k||=\mathcal{O} \left( \sqrt{N}\tau_k \cdot \frac{\alpha_k}{1 - \lambda} \right).\label{17}
\end{equation}
\begin{proof}
See Appendix \ref{sectionD}.\label{lemma7}
\end{proof}
\end{lemma}
\subsection{Convergence Analysis and Main Results}
\hspace{1em}Based on the lemmas shown above, we can derive the theorems bellow.\\
\begin{theorem}\label{theorem1}
Suppose that Assumptions \ref{assumption1}-\ref{assumption8} hold. Let $0<\mu<\frac{1}{2(\rho+1)}$, and assume that the step size $\alpha_k$ and the clipping parameter $\tau_k$ satisfy the following conditions
\begin{equation}
    \begin{split}
      &\sum_{k=0}^{+\infty}\alpha_k=\infty, \lim_{k\to\infty}\alpha_k=0,\lim_{k\to\infty}\frac{\alpha_{k+1}}{\alpha_k}=1\\
      &\sum_{k=0}^{+\infty}\alpha_k^2\tau_k^2<+\infty,\sum_{k=0}^{+\infty}\alpha_k\tau_k^{2-2\alpha}<+\infty\\
      &\tau_k \text{ is increasing and } \lim_{k \to\infty} \tau_k=+\infty.\label{18}
    \end{split}
\end{equation}
Then, $\varphi_{\mu}(\bar{x}_k)$ converges with probability 1.\par
\begin{proof}
See Appendix \ref{setcionE}.
\end{proof}
\end{theorem}
\hspace{1em}Theorem \ref{theorem1} shows that $\varphi_{\mu}(\bar{x}_k)$ converges.
In the following main result, we show that the gradient of $\varphi_{\mu}$ decreases to zero along the sequence $\{\bar{x}_k\}$ under condition \eqref{18}.
\par
\begin{theorem}\label{theorem2}
Under Assumptions \ref{assumption1}-\ref{assumption8}, and supposing that the condition in \eqref{18} holds, we have
\begin{equation}
    \begin{split}
    &\mathbb{E}[\inf_{k=k_0,\dots\infty}||\nabla\varphi_{\mu}(\bar{x}_k)||^2]\\
    &\le \frac{2}{1-2\mu(\rho+1)}\bigg(\frac{||V_{k_0}-\theta||+\sum_{k=k_0}^{\infty}b_k}{\sum_{k=k_0}^{\infty}\alpha_k}\nonumber\\
        &\quad +\frac{\sum_{k=k_0}^{\infty}\frac{ \tau_k^2 \alpha_k^2}{2\mu}+\sum_{k=k_0}^{\infty} \frac{\alpha_k(2\gamma)^{2\alpha}\tau_k^{2-2\alpha}}{2\mu}}{\sum_{k=k_0}^{\infty}\alpha_k}\bigg), \label{27}
    \end{split}
\end{equation}
where \( k_0 \) is the index after which the increasing sequence \( \tau_k \) satisfies \( \tau_k \ge 2C_0 \).\par
\begin{proof}
See Appendix \ref{sectionF}.
\end{proof}
\end{theorem}
\begin{remark}
An example of sequences satisfying \eqref{18} is $\alpha_k=1/(k+1)$ and $\tau_k=2C_0(k+1)^{0.4}$. So when we take the parameters satisfying condition \eqref{18}, we could find an upper bound $M$ such that
\[
    ||V_{k_0}-\theta^*||+\sum_{k=k_0}^{\infty}b_k+\sum_{k=k_0}^{\infty}\frac{\tau_k^2 \alpha_k^2}{2\mu}+\sum_{k=k_0}^{\infty}\frac{\alpha_k(2\gamma)^{2\alpha}\tau_k^{2-2\alpha}}{2\mu}\le M.
\]
Theorem 2 provides the convergence rate 
of $\mathbb{E}[\inf_{k}||\nabla\varphi_{\mu}$\\$(\bar{x}_k)||^2]$. In the above case, for sufficiently large $T$, we have
\[
    \mathbb{E}[\inf_{k_0\le k\le T}||\nabla\varphi_{\mu}(\bar{x}_k)||^2]=\mathcal{O} \left(  \frac{M}{\log{T}}\right).
\]
\end{remark}
\section{NUMERICAL EXPERIMENT}
\hspace{1em}In this section, we apply the proposed algorithm to the phase retrieval problem
\begin{equation}
\min_{\theta\in \Omega}\frac{1}{N}\sum_{i=1}^{N}\bigg(\frac{1}{m}\sum_{j=1}^{m}|\left \langle w_{i,j},\theta \right \rangle^2-y_{i,j} |\bigg). \label{25}
\end{equation}
Following the work of \cite{duchi2019solving}, the goal is to recover the signal $\theta$ under Gaussian measurements $w_{i,j}$, where $w_{i,j}$ denotes the $j$-th row of measurement matrix $W_i$ of agent $i$. The constraint set is given by \( \Omega = \{ \theta \in \mathbb{R}^n \mid \| \theta_i \| \le 1, \ i = 1, \dots, N \} \),  
where \( \theta_i \) represents the local variable associated with agent \( i \).\par
\hspace{1em}Each agent $i$ receives observation $y_{i,j}$ of the form
\[
y_{i,j}=\left\langle w_{i,j},\hat{x} \right\rangle^2,
\]
where $\hat{x}$ denotes the original signal. We rearrange the form as
\begin{equation}
\min_{\theta} \frac{1}{N}\sum_{i=1}^{N}f_i(\theta)=\frac{1}{N}\sum_{i=1}^{N}\frac{1}{m}||(W_i\theta)^2-y_i||_{1}.\label{28}
\end{equation}
The function \(f_i\) has the structure of a composite function \(h(c(\theta))\), where \(h := \|\cdot\|_1\), representing the 1-norm, is a convex function and \(c := (W_i\theta)^2 - y_i\) is a smooth mapping. Therefore, \(f_i\) is a weakly convex function \cite{davis2019stochastic}. In contrast to the works in \cite{chen2021distributed, duchi2019solving, davis2020nonsmooth},  
which focus on deriving convergence rates under the assumptions of sharpness,  
our work investigates the convergence behavior of the stochastic subgradient method with gradient clipping, starting from random initialization.
 Random initialization is model-independent and typically more robust when dealing with model mismatches\cite{chen2019gradient}.\par
\hspace{1em}As suggested by \cite{duchi2019solving}, when $N\times m \ge 3n $ for the proximal linear algorithm, the phase recovery rate is 100$\%$. In our
experiment, we let $N\times m \ge 3n $ for our stochastic subgradient method. The initialization is randomized
and we set $x_{1,0}=x_{2,0}=\dots =x_{N,0}$.\par
\hspace{1em}\textbf{Communication Graph:} To satisfy Assumptions \ref{assumption1}-\ref{assumption3}, we define the elements of adjacency matrix $a_{i,j}$ as follows
\[
a_{ij} =
\begin{cases} 
\frac{1}{3}, & \text{if } j \equiv i , \\
\frac{1}{3}, & \text{if } j \equiv i + 1 , \\
\frac{1}{3}, & \text{if } j \equiv i - 1 , \\
0, & \text{otherwise}.
\end{cases}
\]
This matrix corresponds to a ring structure with 28 nodes, where each node is connected to itself and its two immediate neighbors, with connection weights of $\frac{1}{3}$. Finally, we assume that the adjacency matrix is time-invariant, \emph{i.e.,} $A(k)=A$, the element $\{a_{i,j}\}$ is defined above.\par
\hspace{1em}\textbf{Real-World Experiment:} We use the MNIST dataset, which consists of handwritten digit images \cite{lecun1998gradient}. The signal \( \theta \) in \eqref{25} is recovered using the proposed algorithm with parameters chosen to satisfy condition \eqref{18}.
The image size is $28\times 28$, so we set $n=28\times28, N=28, m=84$, such that the total number of measurements satisfies \( N \times m \ge 3n \).
In other words, each agent owns 84 Gaussian measurements, we take mini-batch $b_{i,k}$ to estimate the subgradient
\begin{align*}
    &g_i(v_{i,k})\\
    &=\frac{1}{b_{i,k}}\sum_{j=1}^{b_{i,k}}\left( 2\langle w_{i,j},v_{i,k} \rangle \cdot \mathrm{sign}(\langle w_{i,j}, v_{i,k} \rangle^2 - y_{i,j}) w_{i,j} \right) .
\end{align*}
In our estimate of subgradient, we investigate the distribution of gradient noise $dist(g_i,\partial f_i(v_{i,k}))$ for each $i$\\
\begin{figure}[htbp]
    \centering
    % 第一张图片
    \begin{subfigure}{0.15\textwidth} % 每张图片占页面宽度的 24%
        \centering
        \includegraphics[width=\textwidth]{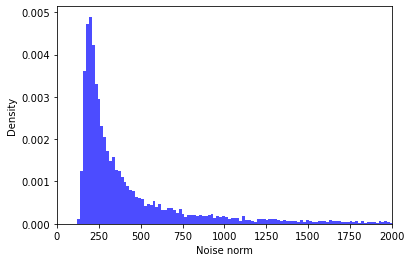} % 替换为你的图片路径
        \caption{Phase retrieval}
    \end{subfigure}
    % 第二张图片
    \begin{subfigure}{0.15\textwidth}
        \centering
        \includegraphics[width=\textwidth]{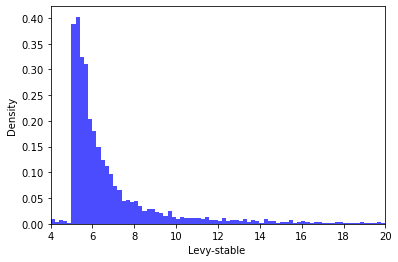}
        \caption{ Lévy-stable}
    \end{subfigure}
    % 第三张图片
    \begin{subfigure}{0.15\textwidth}
        \centering
        \includegraphics[width=\textwidth]{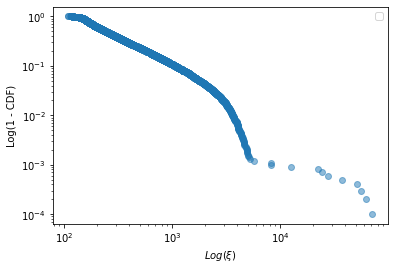}
        \caption{Log-log plot}
    \end{subfigure}

    \caption{(a) Histogram of gradient noise samples for Phase-retrieval on MINIST dataset. (b) Histogram of samples from a sum of squared Lévy-$\alpha$-stable random variables. (c) Log-log plot for subgradient noise.}
    \label{fig:four_images1}
\end{figure}

\hspace{1em}In Fig.\ref{fig:four_images1}, for a fixed agent $i$, we observe that the distribution of the estimate noise  appears heavy-tailed. For comparison, we plot the distribution of Synthetic Lévy-$\alpha$-stable in Fig.\ref{fig:four_images1}(b). Finally, we provide a Log-log plot\cite{gomes1997heavy} to support our observation in Fig.\ref{fig:four_images1}(c).
These results further confirm the heavy-tail behavior of the gradient noise. \par

\begin{figure}[htbp]
    \centering
    % 第一张图片
    \begin{subfigure}{0.1\textwidth} % 每张图片占页面宽度的 24%
        \centering
        \includegraphics[width=\textwidth]{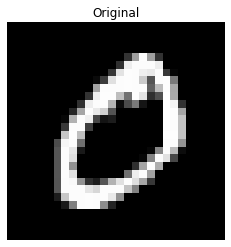} % 替换为你的图片路径
        \caption{}
    \end{subfigure}
    % 第二张图片
    \begin{subfigure}{0.1\textwidth}
        \centering
        \includegraphics[width=\textwidth]{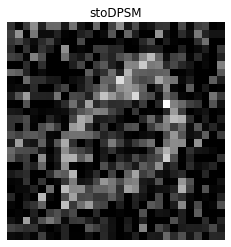}
        \caption{}
    \end{subfigure}
    % 第三张图片
    \begin{subfigure}{0.1\textwidth}
        \centering
        \includegraphics[width=\textwidth]{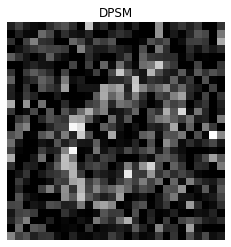}
        \caption{}
    \end{subfigure}
    % 第四张图片
    \begin{subfigure}{0.1\textwidth}
        \centering
        \includegraphics[width=\textwidth]{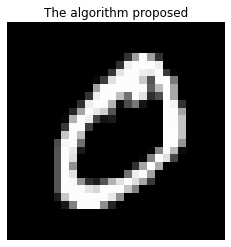}
        \caption{}
    \end{subfigure}

    % 主标题
    \caption{Phase retrieval results on a digit image from the MNIST dataset. From left to right: (a) The original image, (b) The image reconstructed using stoDPSM, (c) The image reconstructed using DPSM, and (d) The image reconstructed using the proposed algorithm. Each method uses $\alpha_k=30/NK$. Data size: $n=784,m=84,N=28.$}
    \label{fig:four_images2}
\end{figure}
\begin{figure}[htbp]
    \centering
    % 第一张图片
    \begin{subfigure}{0.22\textwidth} % 每张图片宽度为 1/4 页面宽度
        \centering
        \includegraphics[width=\textwidth]{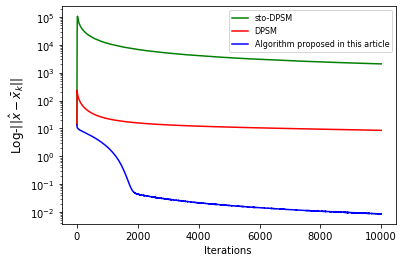} % 替换为你的图片路径
        \caption{}
    \end{subfigure}
    % 第二张图片
    \begin{subfigure}{0.22\textwidth}
        \centering
        \includegraphics[width=\textwidth]{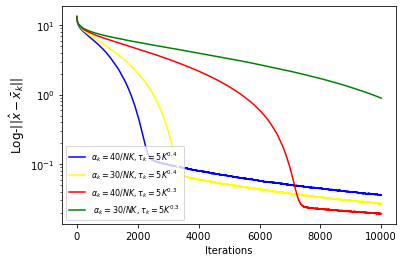}
        \caption{}
    \end{subfigure}
    % 第三张图片
    
    \caption{(a) Linear rate of stoDPSM, DPSM and the proposed algorithm. (b) Linear rate with different parameters.}
    \label{fig:four_images_row}
\end{figure}

\hspace{1em}Motivated by the observed heavy-tail behavior of the gradient noise,  
we compare the proposed algorithm with two baselines: DPSM and stoDPSM, both implemented without spectral initialization\cite{lecun1998gradient}. In Fig.\ref{fig:four_images2}, we compare the original image with those generated by different methods. After the recovery process, our method aligns closely with the original one. Fig. \ref{fig:four_images_row} illustrates the convergence curves of DPSM, stoDPSM, and the proposed algorithm.

\section{CONCLUSION}
\hspace{1em}In our work, we study the distributed stochastic weakly convex optimization under heavy-tailed noise. We propose an update law using gradient clipping, and utilize Moreau envelope
as a potential function to prove convergence. Unlike previous works on distributed stochastic optimization, this paper assumes that gradient noise follows a heavy-tailed distribution, an  extended and more general assumption than the commonly used bounded variance noise.
Moreover, we assume that the objective function is non-smooth and weakly convex. By numerical experiments on robust phase retrieval, we investigate the heavy-tail phenomenon and demonstrate the effectiveness of our algorithm.

\section{APPENDIX}
\appendix
\section{Proof of \eqref{6}}\label{AppendixA}
\hspace{1em}The proof is motivated by \cite{boyd2013proximal}. Recalling that the Moreau envelope is given by
\[
\varphi_{\mu}(x)=\min_y{\varphi(y)+\frac{1}{2\mu}||y-x||^2},
\]
we rearrange the above term as follows
\begin{align*}
    \varphi_{\mu}(x)&=\frac{1}{2\mu}||x||^2-\frac{1}{\mu}\sup_{y}\{x^{\top}y-\mu f(y)-\frac{1}{2}||y||^2\}\\
                &=\frac{1}{2\mu}||x||^2-\frac{1}{\mu}\bigg(\mu f+\frac{1}{2}||\cdot||^2\bigg)^*(x),
\end{align*}
where the second equation holds by the definition of conjugate function\cite{klee1971convex}. Since in the compact set case, we have $inf*=min*$, so we change the expression.\par
\hspace{1em} By \cite{boyd2013proximal}, we have that the gradient of the conjugate function
is equal to the optimal $y^*$ as follows
\begin{align*}
    \nabla\varphi_{\mu}(x)&=\frac{x}{\mu}-\frac{1}{\mu}\mathop{\text{argmax}}\limits_{y}\{x^{\top}y-\mu f(y)-\frac{1}{2}||y||^2\}\\
                      &=\frac{1}{\mu}(x-prox_{\mu f}(x))=\frac{1}{\mu}(x-\hat{x}).
\end{align*}

\section{Proof of Lemma 3}\label{sectionB}
\hspace{1em}Taking the norm of the bias, and according to Jensen's inequality, we have
    \[
    \Vert \mathbb{E}[\hat{g}_{i,k}|\mathcal{F}_k]-G_{i,k}\Vert\le\mathbb{E}[\vert \Vert g_{i,k} \Vert -\tau_k \vert \mathbf{1}_{\{||g_{i,k}||>\tau_k\}}|\mathcal{F}_k].
    \]
    By the triangle inequality, we further derive that
    \begin{align*}
    &\mathbb{E}[ \vert \Vert g_{i,k} \Vert -\tau_k \vert\mathbf{1}_{\{||g_{i,k}||>\tau_k\}}|\mathcal{F}_k]\\
    &\le\mathbb{E}[ (\Vert  g_{i,k}\Vert  -C_0) \mathbf{1}_{\{||g_{i,k}||>\tau_k\}}|\mathcal{F}_k]\\
    &  \le\mathbb{E}[||g_{i,k}-G_{i,k}||\mathbf{1}_{\{||g_{i,k}||>\tau_k\}}|\mathcal{F}_k].
    \end{align*}
    Since ${\{||g_{i,k}||>\tau_k\}} \subset \{||g_{i,k}-G_{i,k}||>\frac{\tau_k}{2}\}$, then $\mathbf{1}_{\{||g_{i,k}||>\tau_k\}}\le \mathbf{1}_{\{||g_{i,k}-G_{i,k}||>\frac{\tau_k}{2}\}}$. Hence, we have
    \begin{align*}
        &\mathbb{E}[||g_{i,k}-G_{i,k}||\mathbf{1}_{\{||g_{i,k}||>\tau_k\}}|\mathcal{F}_k]\\
        &  \le\mathbb{E}[||g_{i,k}-G_{i,k}||\mathbf{1}_{\{||g_{i,k}-G_{i,k}||>\frac{\tau_k}{2}\}}|\mathcal{F}_k].
    \end{align*}
    Based on the conditional Holder’s inequality,
    \begin{align*}
        ||B_{i,k}||&\le(E[||g_{i,k}-G_{i,k}||^{\alpha}|\mathcal{F}_k])^{\frac{1}{\alpha}}\\
                   &\quad\cdot(E[\mathbf{1}_{\{||g_{i,k}-G_{i,k}||>\frac{\tau_k}{2}\}}^{\frac{\alpha}{\alpha-1}}|\mathcal{F}_k])^{\frac{\alpha-1}{\alpha}}\\
                   &=\gamma(E[\mathbf{1}_{\{||g_{i,k}-G_{i,k}||>\frac{\tau_k}{2}\}}|\mathcal{F}_k])^{\frac{\alpha-1}{\alpha}}.
    \end{align*}
    Applying Markov's inequality to conditional expectation, 
    \begin{align*}
        \mathbb{E}[\mathbf{1}_{\{||g_{i,k}-G_{i,k}||>\frac{\tau_k}{2}\}}|\mathcal{F}_k]&=prob(||g_{i,k}-G_{i,k}||>\frac{\tau_k}{2}|\mathcal{F}_k)\\
                                                                          &=prob(||g_{i,k}-G_{i,k}||^{\alpha}>\frac{\tau_k^{\alpha}}{2^{\alpha}}|\mathcal{F}_k)\\
                                                                          &\le\frac{\mathbb{E}[||g_{i,k}-G_{i,k}||^{\alpha}|\mathcal{F}_k]}{\frac{\tau_k^{\alpha}}{2^{\alpha}}}.
    \end{align*}
    Finally,
    \begin{align*}
    &||B_{i,k}||\le (2\gamma)^\alpha\tau_k^{1-\alpha}.
    \end{align*}

\section{Proof of Lemma 4}\label{sectionC}
\hspace{1em}The following inequality holds due to the nonexpansiveness of the projection operator
    \begin{align*}
        \| x_{i,k+1} - \hat{v}_{i,k} \|^2 &=||\mathcal{P}_{\Omega}[v_{i,k}-\alpha_k\hat{g}_{i,k}]-\hat{v}_{i,k}||^2\\
        &\leq||v_{i,k}-\hat{v}_{i,k}-\alpha_k\hat{g}_{i,k}||^2\\
        &= \| v_{i,k} - \hat{v}_{i,k} \|^2-2\alpha_k\hat{g}_{i,k}^{\top}(v_{i,k}-\hat{v}_{i,k})\\
        &\quad+\alpha_k^2||\hat{g}_{i,k}||^2\\
        &=|| v_{i,k} - \hat{v}_{i,k}||^2-2\alpha_k G_{i,k}^{\top}(v_{i,k}-\hat{v}_{i,k})\\
        &\quad-2\alpha_kb_{i,k}^{\top}(v_{i,k}-\hat{v}_{i,k})+\alpha_k^2 ||\hat{g}_{i,k}||^2.
        \end{align*}
Recalling the weak convexity of $f_i$ in \eqref{1} and the boundness of $\hat{g}_{i,k}$, it gives that
    \begin{align*}
        ||x_{i,k+1}-\hat{v}_{i,k}||^2&\le||v_{i,k}-\hat{v}_{i,k}||^2+2\alpha_k(f_i(\hat{v}_{i,k})-f_i(v_{i,k})\\
                                     &\quad+\frac{\rho}{2}||v_{i,k}-\hat{v}_{i,k}||^2)-2\alpha_kb_{i,k}^{\top}(v_{i,k}-\hat{v}_{i,k}) \\
                                     &\quad+\alpha_k^2\tau_k^2.
    \end{align*}
    Then, we consider the upper bound of each term
    \begin{equation}
        \begin{split}
           &f_i(\hat{v}_{i,k}) - f_i(v_{i,k})\\
           &= f_i(\hat{v}_{i,k}) - f_i(z_k) + f_i(z_k) - f_i(\bar{x}_k) + f_i(\bar{x}_k) - f_i(v_{i,k})\\
           &\leq L \left( \frac{1}{1 - \mu\rho} + 1 \right) \| v_{i,k} - \bar{x}_k \| + f_i(z_k) - f_i(\bar{x}_k)\\
           &\leq \frac{L(2 - \mu\rho)}{1 - \mu\rho} \sum_{j=1}^{N} a_{i,j}(k)\| x_{j,k} - \bar{x}_k \| + f_i(z_k) - f_i(\bar{x}_k),\label{26}
        \end{split}
    \end{equation}
    where the first inequality holds due to Lemma \ref{lemma2} and assumption \ref{assumption6}, and the last inequality holds because of the convexity of the norm. Using the inequality $||a+b||^2\le2||a||^2+2||b||^2$,
    \begin{equation}
        \begin{split}
          & \|v_{i,k} - \hat{v}_{i,k}\|^2 \\
          &=  \|v_{i,k} - \bar{x}_k + \bar{x}_k - z_k + z_k - \hat{v}_{i,k}\|^2 \\
          &\leq  2\|\bar{x}_k - z_k\|^2 +  2\|v_{i,k} - \bar{x}_k + z_k - \hat{v}_{i,k}\|^2 \\
          &\leq  2\|\bar{x}_k - z_k\|^2 +  4\left( 1 + \frac{1}{(1 - \mu\rho)^2} \right) \|v_{i,k} - \bar{x}_k\|^2 \\
          &\leq  2\|\bar{x}_k - z_k\|^2 \\
          &\quad+ 4\left( 1 + \frac{1}{(1 - \mu\rho)^2} \right) \sum_{j=1}^{N} a_{i,j}(k)\|x_{j,k} - \bar{x}_k\|^2.\label{27}
        \end{split}  
    \end{equation}
        According to Assumption \ref{assumption1}, summing  \eqref{26} and $\rho/2$ times \eqref{27} from $i=1$ to $N$ yields
        \begin{align*}
            &\sum_{i=1}^{N} \left( f_i(\hat{v}_{i,k}) - f_i(v_{i,k}) + \frac{\rho}{2} \|v_{i,k} - \hat{v}_{i,k}\|^2 \right) \\
            &\leq \frac{L(2 - \mu\rho)}{1 - \mu\rho} \sum_{i=1}^{N} \|x_{i,k} - \bar{x}_k\| + N(f(z_k) - f(\bar{x}_k)) \\
            &\quad + N \rho \|\bar{x}_k - z_k\|^2 + 2 \rho \left( 1 + \frac{1}{(1 - \mu\rho)^2} \right) \sum_{i=1}^{N} \|x_{i,k} - \bar{x}_k\|^2.
        \end{align*}
        From the definition of $z_k$, if $t<\frac{1}{2\rho}$, one has
        \begin{align*}
            &f(z_k)-f(\bar{x}_k)+\rho||\bar{x}_k-z_k||^2\\
            &=f(z_k)-f(\bar{x}_k)+(\frac{1}{2\mu}-\frac{1}{2\mu}+\rho)||\bar{x}_k-z_k||^2\\
            &\leq(-\frac{1}{2\mu}+\rho)||\bar{x}_k-z_k||^2.
        \end{align*}  
        Finally, we obtain
        \begin{align*}
        &\sum_{i=1}^{N}||x_{i,k+1}-\hat{v}_{i,k}||\\
        \leq&\sum_{i=1}^{N}||v_{i,k}-\hat{v}_{i,k}||^2\\
            &\ \ +2\alpha_k (N\left(-\frac{1}{2\mu}+\rho\right)||\bar{x}_k-z_k||^2\\
            &\ \ +\frac{L(2-\mu\rho)}{1-\mu\rho}\sum_{i=1}^{N}||\bar{x}_k-x_{i,k}||\\
            &\ \ +2\rho\left(1+\frac{1}{(1-\mu\rho)^2}\right)\sum_{i=1}^{N}||\bar{x}_k-x_{i,k}||^2\\
            &\ \ -\sum_{i=1}^{N}b_{i,k}^{\top}(v_{i,k}-\hat{v}_{i,k}))+N\tau_k^2\alpha_k^2.
        \end{align*} 

\section{Proof of Lemma 5}\label{sectionD}

    \hspace{1em}Motivated by \cite{liu2017convergence}, we assume that \( n=1 \) and define
    \begin{align*}
        x_k &= [x_{1,k},x_{2,k},\dots,x_{N,k}]^{\top}\\
        v_k &= [v_{1,k},v_{2,k},\dots,v_{N,k}]^{\top}\\
        \epsilon_k &= [e_{1,k},e_{2,k},\dots,e_{N,k}]^{\top},
    \end{align*}
    where \( \epsilon_{i,k} = \mathcal{P}_{\Omega}[v_{i,k}-\alpha_k\hat{g}_{i,k}]-v_{i,k} \). Rewriting the iteration \eqref{14}, we obtain
    \[
    x_{k+1} = v_k + \epsilon_k = A(k)x_k + \epsilon_k.
    \]
    From the definition of \( \hat{g}_{i,k} \), it follows that
    \[
    \|\epsilon_{i,k}\|^2 \leq \|v_{i,k}- v_{i,k} - \alpha_k \hat{g}_{i,k}\|^2 \leq \alpha_k^2 \tau_k^2.
    \]
    Thus, summing over all agents, we obtain that
    \[
    \|\epsilon_k\| \leq \sqrt{N} \alpha_k \tau_k.
    \]
    Applying Lemma 3 and following the proof of \cite[Lemma II.5]{nedic2010constrained}, for some constants \( c>0 \) and \( \lambda \in (0,1) \), independent of \( k \), it holds that
    \[
    \|\Delta_{k+1}\| \leq c\lambda^k \|\Delta_0\| + c\sqrt{N} \tau_k \sum_{l=0}^{k-1} \lambda^{k-l-1} \alpha_l + \sqrt{N} \alpha_k \tau_k,
    \]
    where the $i$-th element of $\Delta_{k}$ equals to $||x_{i,k}-\bar{x}_k||$.
    Finally, applying Lemma 4 and the limit condition \( \lim_{k\to\infty} \alpha_{k+1}/\alpha_k = 1 \), the result follows.
    
\section{Proof of Theorem 1}\label{setcionE}
\hspace{1em}By the definition of $\varphi_{\mu}$ in \eqref{3}, we obtain 
\[
\varphi_{\mu}(x_{i,k+1})=\mathop{\text{min}}\limits_{y \in \Omega}f(y)+\frac{1}{2\mu}||x_{i,k+1}-y||^2.
\]
Let $\hat{v}_{i,k}=\mathop{\text{argmin}}\limits_{y \in \Omega}f(y)+\frac{1}{2\mu}||y-v_{i,k}||^2$, $\hat{x}_{i,k}=\mathop{\text{argmin}}\limits_{y \in \Omega}f(y)+\frac{1}{2\mu}||y-x_{i,k}||^2$. In the case of $y=\hat{v}_{i,k}$, we obtain that
\begin{equation}
\varphi_{\mu}(x_{i,k+1})\le f(\hat{v}_{i,k})+\frac{1}{2\mu}||x_{i,k+1}-\hat{v}_{i,k}||^2.\label{19}
\end{equation}
Summing \eqref{19} from $1$ to $N$ and  using  Lemma \ref{lemma6} yields
\begin{equation}\label{20}
    \begin{split}
      \sum_{i=1}^N \varphi_{\mu} (x_{i,k+1})  
      &\leq \sum_{i=1}^N \left(f (\hat{v}_{i,k})+ \frac{1}{2\mu}||v_{i,k}-\hat{v}_{i,k}||^2\right)\\
      &\quad + \frac{\alpha_k}{\mu}  \bigg(N\left(-\frac{1}{2\mu} + \rho\right) \| \bar{x}_k - z_k \|^2\\
      &\quad + \frac{L(2 - \mu\rho)}{1 - \mu\rho} \sum_{i=1}^N \| x_{i,k} - \bar{x}_k \| \\
      &\quad + 2\rho \bigg( 1 + \frac{1}{(1 - \mu\rho)^2} \bigg)\sum_{i=1}^N \| x_{i,k} - \bar{x}_k \|^2 \\
      &\quad -\sum_{i=1}^{N}b_{i,k}^{\top}(v_{i,k}-\hat{v}_{i,k})\bigg)+\frac{N \tau_k^2 \alpha_k^2}{2\mu}.  
    \end{split} 
\end{equation}
\hspace{1em}Since  $\varphi_{\mu}(v_{i,k})= \min_{y\in \Omega}f(y)+\frac{1}{2\mu}||v_{i,k}-y||^2$ = $f(\hat{v}_{i,k})+\frac{1}{2\mu}||v_{i,k}-\hat{v}_{i,k}||^2$ and
noting that $v_{i,k}=\sum_{j=1}^{N}a_{i,j}(k)x_{j,k}$, we have
\begin{equation}
    \begin{split}
    &\varphi_{\mu}(v_{i,k})\\
    &\quad=f\left(\sum_{j=1}^{N}a_{i,j}(k)\hat{v}_{i,k}\right)+\frac{1}{2\mu}||\sum_{j=1}^{N}a_{i,j}(k)(\hat{v}_{i,k}-x_{j,k})||^2\\
    &\quad\le f\left(\sum_{j=1}^{N}a_{i,j}(k)\hat{x}_{j,k}\right)+\frac{1}{2\mu}\sum_{j=1}^{N}a_{i,j}(k)||\hat{x}_{j,k}-x_{j,k}||^2\\
    &\quad\le \sum_{j=1}^{N}a_{i,j}(k)f(\hat{x}_{j,k})\\
    &\qquad +\frac{\rho}{2}\sum_{1\le j \le l \le N}a_{i,j}(k)a_{i,l}(k)||\hat{x}_{j,k}-\hat{x}_{l,k}||^2\\
    &\qquad +\frac{1}{2\mu}\sum_{j=1}^{N}a_{i,j}(k)||\hat{x}_{j,k}-x_{j,k}||^2\\
    &\quad \le\sum_{j=1}^{N}a_{i,j}(k)\varphi_{\mu}(x_{j,k})\\
    &\qquad +\frac{\rho}{2(1-\mu\rho)^2}\sum_{1\le j \le l \le N}a_{i,j}(k)a_{i,l}(k)||{x}_{j,k}-{x}_{l,k}||^2,\\
    \end{split}\label{70}
\end{equation}
where the first inequality holds by the definition of $\varphi_{\mu}(\hat{v}_{i,k})$, the second one follows Lemma \ref{lemma1} and the last inequality holds
because of Lemma \ref{lemma2}. Letting $V_{k+1}=1/N\sum_{i=1}^{N}\varphi_{\mu}(x_{i,k+1})$ and following from \eqref{20} and \eqref{70}, we have
\begin{equation}
  \begin{split}
    &V_{k+1}\le V_{k}\\
    &\quad+\frac{\rho}{2N(1-\mu\rho)^2}\cdot\\
    &\qquad \sum_{i=1}^{N}\sum_{1\le j \le l \le N}a_{i,j}(k)a_{i,l}(k)||{x}_{j,k}-{x}_{l,k}||^2\\
    &\quad + \frac{\alpha_k}{\mu}  \bigg(\bigg(\rho-\frac{1}{2\mu} \bigg) \| \bar{x}_k - z_k \|^2\\
    &\quad + \frac{L(2 - \mu\rho)}{N(1 - \mu\rho)} \sum_{i=1}^N \| x_{i,k} - \bar{x}_k \| \\
    &\quad + \frac{2\rho}{N} \bigg( 1 + \frac{1}{(1 - \mu\rho)^2} \bigg)\sum_{i=1}^N \| x_{i,k} - \bar{x}_k \|^2 \\
    &\quad -\frac{1}{N}\sum_{i=1}^{N}b_{i,k}^{\top}(v_{i,k}-\hat{v}_{i,k})\bigg)+\frac{ \tau_k^2 \alpha_k^2}{2\mu}.
  \end{split}\label{71}
\end{equation}
\hspace{1em}Taking the $\sigma-$algebra of inequality \eqref{71} yields
\begin{equation}
    \begin{split}
      &\mathbb{E}[V_{k+1}|\mathcal{F}_k]\le V_{k}\\
      &\quad+\frac{\rho}{2N(1-\mu\rho)^2}\cdot\\
      &\qquad \sum_{i=1}^{N}\sum_{1\le j \le l \le N}a_{i,j}(k)a_{i,l}(k)||{x}_{j,k}-{x}_{l,k}||^2\\
      &\quad + \frac{\alpha_k}{\mu}  \bigg(\bigg(\rho-\frac{1}{2\mu}\bigg) \| \bar{x}_k - z_k \|^2\\
      &\quad + \frac{L(2 - \mu\rho)}{N(1 - \mu\rho)} \sum_{i=1}^N \| x_{i,k} - \bar{x}_k \| \\
      &\quad + \frac{2\rho}{N} \bigg( 1 + \frac{1}{(1 - \mu\rho)^2} \bigg)\sum_{i=1}^N \| x_{i,k} - \bar{x}_k \|^2 \\
      &\quad -\frac{1}{N}\sum_{i=1}^{N}B_{i,k}^{\top}(v_{i,k}-\hat{v}_{i,k})\bigg)+\frac{ \tau_k^2 \alpha_k^2}{2\mu}. \label{21}
    \end{split}
\end{equation}
\hspace{1em}Since $\tau_k$ is increasing and $\lim_{k\to\infty} \tau_k=+\infty$, there exists $k_o$ such that for $k \ge k_0, \tau_k\ge 2C_0 $. Thus, Lemma \ref{lemma5} holds for $k\ge k_0$. Applying Lemma \ref{lemma5} and Young's inequality to the term $\frac{1}{N}\sum_{i=1}^{N}B_{i,k}^{\top}(v_{i,k}-\hat{v}_{i,k})$ gives
\begin{equation}
    \begin{split}
      -\frac{1}{N}\sum_{i=1}^{N}B_{i,k}^{\top}&(v_{i,k}-\hat{v}_{i,k})\\
                                                   &\le\frac{1}{2N}\sum_{i=1}^{N}\left(||B_{i,k}||^2+||v_{i,k}-\hat{v}_{i,k}||^2\right)\\
                                                   &\le\frac{1}{2}(2\gamma)^{2\alpha}\tau_k^{2-2\alpha}+\frac{1}{2N}\sum_{i=1}^{N}||v_{i,k}-\hat{v}_{i,k}||^2\\
                                                   &\le\frac{1}{2}(2\gamma)^{2\alpha}\tau_k^{2-2\alpha}+ ||\bar{x}_k-z_k||^2\\
                                                   &\quad +\frac{2}{N} \left( 1 + \frac{1}{(1 - \mu\rho)^2} \right) \sum_{i=1}^{N} \|x_{i,k} - \bar{x}_k\|^2,\label{22}
    \end{split}
\end{equation}
where the last inequality holds based on \eqref{27} in Appendix \ref{sectionB}.\par

\hspace{1em} Simplifying the right side of \eqref{21} gives
\[
    \mathbb{E}[V_{k+1}|\mathcal{F}_k]\le V_{k}+b_k+\frac{\alpha_k(2\gamma)^{2\alpha}\tau_k^{2-2\alpha}}{2\mu}+\frac{ \tau_k^2 \alpha_k^2}{2\mu}
\]
where 
\begin{align*}
     b_k:=&\quad A_{\mu\rho} 
    \sum_{i=1}^{N}\sum_{1\le j\le l \le N} a_{i,j}(k)a_{i,l}(k)\| {x}_{j,k} - {x}_{l,k} \|^2 \nonumber \\
    &\quad + \frac{\alpha_k}{\mu} \bigg( B_{\mu\rho}\sum_{i=1}^N \| x_{i,k} - \bar{x}_k \|  \nonumber \\
    &\quad + C_{\mu\rho} \sum_{i=1}^N \| x_{i,k} - \bar{x}_k \|^2 \bigg)
\end{align*}
and $A_{\mu\rho}= \frac{\rho}{2N(1-\mu\rho)^2}$, $B_{\mu\rho}:=\frac{L(2 - \mu\rho)}{N(1 - \mu\rho)} $, $C_{\mu\rho}=\frac{2\rho+2}{N} \bigg(1 + \frac{1}{(1 - \mu\rho)^2}\bigg)$ and we used the condition  $-\frac{1}{2\mu}+\rho+1<0$.\par
\hspace{1em}To analyze the convergence, we need to derive the upper bound of $b_k$.  Based on Lemma \ref{lemma7},
\begin{align*}
    \sum_{i=1}^{N}&\sum_{1\le j \le l \le N} a_{i,j}(k)a_{i,l}(k)\| {x}_{j,k} - {x}_{l,k} \|^2= \mathcal{O} \left(  \frac{2N\tau_k^2 \alpha_k^2}{(1 - \lambda)^2} \right)\\
    &\alpha_k\sum_{i=1}^{N}\| x_{i,k} - \bar{x}_k \|^2=\mathcal{O} \left(  \frac{N\tau_k^2 \alpha_k^3}{1 - \lambda} \right)\\
    &\alpha_k\sum_{i=1}^{N}\| x_{i,k} - \bar{x}_k \|=\mathcal{O} \left(  \frac{N\tau_k \alpha_k^2}{1 - \lambda} \right),
\end{align*}
where the first and the second equation hold  since $||\Delta_k||^2=\sum_{i=1}^{N}||x_{i,k}-\bar{x}_k||^2=\mathcal{O} \left( N\tau_k^2 \cdot \frac{\alpha_k^2}{(1 - \lambda)^2} \right)$, and the last equation holds because of Cauchy-Schwarz inequality,   $\bigg(\sum_{i=1}^N ||x_{i,k}-\bar{x}_k||\bigg)^2\le N\sum_{i=1}^N ||x_{i,k}-\bar{x}_k||^2$.\par
\hspace{1em}Thus, $b_k=\mathcal{O} \left(  2N\frac{\tau_k^2 \alpha_k^2}{(1 - \lambda)^2} \right)$. By \cite[Theorem 3.1]{dhasunr2010distributed}, $V_{k}$ converges to a 
nonegative random variable $\theta$ with probability 1.

\hspace{1em}The definition of Moreau envelope in \eqref{3} states that $\varphi_{\mu}(x) $ is $C_1$ smooth. Hence, by applying Lemma \ref{lemma7} with $||\bar{x}_k-x_{i,k}||\to0$, it follows that
\begin{equation}
|\varphi_{\mu}(x_{i,k})-\varphi_{\mu}(\bar{x}_k)|\to0, \label{23}
\end{equation}
and\\
\begin{equation}
    \begin{split}
    \vert V_{k}-\varphi_{\mu}(\bar{x}_k) \vert=& \vert \frac{1}{N}\sum_{i=1}^N\varphi_{\mu}(x_{i,k})-\varphi_{\mu}(\bar{x}_k) \vert\\
                                            \le&\frac{1}{N}\sum_{i=1}^N|\varphi_{\mu}(x_{i,k})-\varphi_{\mu}(\bar{x}_k)\vert \to0. \label{24}
    \end{split}
\end{equation}
Thus, $lim_{k\to\infty}V_{k}=lim_{k\to\infty}\varphi_{\mu}(\bar{x}_k)$. Moreau envelope $\varphi_{\mu}(\bar{x}_k)$ also converges.\par
\section{Proof of Theorem 2}\label{sectionF}
\hspace{1em}The inequality \eqref{21} can be rewritten as
\begin{equation}
    \begin{split}
    \frac{\alpha_k}{\mu}\left(\frac{1}{2\mu}-\rho-1\right)||\bar{x}_k-z_k||^2&\le V_{k}-\mathbb{E}[V_{k+1}|\mathcal{F}_k]\\
                                                                   &+b_k+\frac{\alpha_k(2\gamma)^{2\alpha}\tau_k^{2-2\alpha}}{2\mu}\\
                                                                   &+\frac{ \tau_k^2 \alpha_k^2}{2\mu}.
    \end{split}
\end{equation}
Summing up both sides from $k=k_0$ to infinity gives
\begin{align*}
    \frac{1-2\mu(\rho+1)}{2\mu^2}\sum_{k=k_0}^{\infty}\alpha_k||\bar{x}_k-z_k||^2\le&\sum_{k=k_0}^{\infty}(V_{k}-\mathbb{E}[V_{k+1}|\mathcal{F}_k])\\
                                                     &+\sum_{k=k_0}^{\infty}b_k+\sum_{k=k_0}^{\infty}\frac{ \tau_k^2 \alpha_k^2}{2\mu}\\
                                                     &+ \sum_{k=k_0}^{\infty} \frac{\alpha_k(2\gamma)^{2\alpha}\tau_k^{2-2\alpha}}{2\mu}.                                 
\end{align*}
Based on equation \eqref{6}, we have $||\bar{x}_k-z_k||=\mu \nabla\varphi_{\mu}(\bar{x}_k)$.
Then dividing both sides by $\sum_{k=k_0}^{\infty}\alpha_k$, we obtain
\begin{align*}
    &\inf_{k=k_0,...,\infty}||\nabla\varphi_{\mu}(\bar{x}_k)||^2\\
    &\le \frac{2}{1-2\mu(\rho+1)}\bigg(\frac{\sum_{k=k_0}^{\infty}(V_{k}-\mathbb{E}[V_{k+1}|\mathcal{F}_k])+\sum_{k=k_0}^{\infty}b_k}{\sum_{k=k_0}^{\infty}\alpha_k}\nonumber\\
    &\quad +\frac{\sum_{k=k_0}^{\infty}\frac{ \tau_k^2 \alpha_k^2}{2\mu}+\sum_{k=k_0}^{\infty} \frac{\alpha_k(2\gamma)^{2\alpha}\tau_k^{2-2\alpha}}{2\mu}}{\sum_{k=k_0}^{\infty}\alpha_k}\bigg).\\
\end{align*}
Taking expectation on both sides gives
\begin{align*}
&\mathbb{E}[\inf_{k=k_0,...,\infty}||\nabla\varphi_{\mu}(\bar{x}_k)||^2]\\
&\le \frac{2}{1-2\mu(\rho+1)}\bigg(\frac{||V_{k_0}-\theta||+\sum_{k=k_0}^{\infty}b_k}{\sum_{k=k_0}^{\infty}\alpha_k}\nonumber\\
&\quad +\frac{\sum_{k=k_0}^{\infty}\frac{ \tau_k^2 \alpha_k^2}{2\mu}+\sum_{k=k_0}^{\infty} \frac{\alpha_k(2\gamma)^{2\alpha}\tau_k^{2-2\alpha}}{2\mu}}{\sum_{k=k_0}^{\infty}\alpha_k}\bigg).\\
\end{align*}
\bibliographystyle{unsrt}
\bibliography{reference.bib}
\end{document}